\documentclass[12pt]{article}
\usepackage{latexsym}
\usepackage{amssymb}
\usepackage{amsmath}

\begin{document}

\vspace*{.5cm}
\begin{center}
{\large
{\bf Invertibility Preserving Linear Maps\\ On the Semi-Simple
Banach Algebras}}
\bigskip

{\bf Mohamad Reza Farmani}
\\
{\tt mr.farmanis@gmail.com}
\end{center}
\bigskip

\begin{abstract}
In this paper, we show that the essentiality of the scole of an
ideal $\cal{B}$ in a semi-simple Banach algebra $\cal{A}$ implies
that any invertibility preserving isomorphism $\phi:
{\cal{A}}\rightarrow{\cal{A}}$ is a Jordan homomorphism. Specially
if, the unitary semi-simple Banach algebra $\cal{A}$ has an
essential minimal ideal
then $\phi\mid_{soc({\cal{A}})}$ is a Jordan homomorphism.\\

\end{abstract}


\newpage


\begin{thebibliography}{99}
\bibitem{}  B. Aupetit, {\it  A primer on spectral theory,}
Springer-Verlag, New York, 1991.
\bibitem{} M. Bresar, A. Fosner and P. Semrl,  {\it
A note on invertibility preserving Banach algebra, } American
Math. Soc.,{\bf 131}(2003), 3833-3837.
\bibitem{} A. Harris, {\it Invertibility preserving linear maps on Banach
algebras,} Israel Mathematical Conferance Proceedings, {\bf 15}
(2005), 156-171
\bibitem{} M. Marthieu and C. Ruddy, {\it Spectral isometries},
Contemporary Mathematics, 2006.
\end{thebibliography}
\end{document}